# Some remarks on generalized Fibonacci and Lucas polynomials

*Johann Cigler*

Fakultät für Mathematik, Universität Wien

**Abstract**

Starting with some determinants of binomial coefficients which are related to Fibonacci and Lucas polynomials we study similar determinants for some generalizations of these polynomials and their q-analogues.

## 1. Introduction

The Fibonacci polynomials $F_n(x,s) = \sum_{k=0}^{\lfloor n/2 \rfloor} \binom{n-j}{j} s^j x^{n-2j}$ satisfy the recursion

$F_n(x,s) = x F_{n-1}(x,s) + s F_{n-2}(x,s)$ with initial values $F_0(x,s) = 1$ and $F_1(x,s) = x$.

The Lucas polynomials $L_n(x,s) = \sum_{k=0}^{\lfloor n/2 \rfloor} \frac{n}{n-j} \binom{n-j}{j} s^j x^{n-2j}$ satisfy the same recursion

$L_n(x,s) = x L_{n-1}(x,s) + s L_{n-2}(x,s)$ but with initial values $L_0(x,s) = 2$ and $L_1(x,s) = x$.

For $s=1$ we write $F_n(x,1) = F_n(x)$ and $L_n(x,1) = L_n(x)$.

The starting point of this paper was the observation that the Fibonacci and Lucas polynomials can be obtained as the following determinants of matrices with binomial coefficients:

$$\det\left(\binom{i-1}{j} x^2 + \binom{i+1}{j+1}\right)_{i,j=0}^{n-1} = F_{2n}(x), \tag{1}$$

$$x \det\left(\binom{i}{j} x^2 + \binom{i+2}{j+1}\right)_{i,j=0}^{n-1} = F_{2n+1}(x), \tag{2}$$

$$\det\left(\binom{2i}{i-j} x^2 + \binom{2i+2}{i+1-j}\right)_{i,j=0}^{n-1} = L_{2n}(x), \tag{3}$$

$$x \det\left(\binom{2i+1}{i-j} x^2 + \binom{2i+3}{i+1-j}\right)_{i,j=0}^{n-1} = L_{2n+1}(x), \tag{4}$$

$$\det\left(\binom{i+1}{j+1} x^2 - \binom{i}{j-1}\right)_{i,j=0}^{n-1} = x^n F_n(x). \tag{5}$$

Note that we always assume that $\binom{x}{k} = 0$ for $k < 0$. It is very probable that these determinants are known, but I could find no reference.

We extend these results to the generalized Fibonacci polynomials

$$F_n^{(k)}(x,s) = \sum_{j=0}^{\left\lfloor \frac{n}{k} \right\rfloor} \binom{n-(k-1)j}{j} s^j x^{n-kj} \tag{6}$$

which satisfy the recursion $F_n^{(k)}(x,s) = xF_{n-1}^{(k)}(x,s) + sF_{n-k}^{(k)}(x,s)$ with initial values $F_n^{(k)}(x,s) = x^n$ for $0 \leq n < k$ and to the generalized Lucas polynomials

$$L_n^{(k)}(x,s) = \sum_{j=0}^{\left\lfloor \frac{n}{k} \right\rfloor} \frac{n}{n-(k-1)j} \binom{n-(k-1)j}{j} s^j x^{n-kj} \tag{7}$$

which satisfy the same recursion $L_n^{(k)}(x,s) = xL_{n-1}^{(k)}(x,s) + sL_{n-k}^{(k)}(x,s)$ but with initial values $L_0^{(k)}(x,s) = k$ and $L_n^{(k)}(x,s) = x^n$ for $0 < n < k$.

Note that $F_n^{(1)}(x,s) = L_n^{(1)}(x,s) = (x+s)^n$ and $F_n^{(2)}(x,s) = F_n(x,s)$ and $L_n^{(2)}(x,s) = L_n(x,s)$.

These polynomials are related by

$$L_{n+k}^{(k)}(x,s) = F_{n+k}^{(k)}(x,s) + (k-1)sF_n^{(k)}(x,s) \tag{8}$$

for $n \in \mathbb{N}$ as is easily verified by induction.

The recursion can be used to extend these polynomials to all $n \in \mathbb{Z}$. We shall need only the fact that they can be extended to $-k < n < 0$ with $F_n^{(k)}(x,s) = 0$.

The numbers $F_n^{(k)}(1,1)$ and $L_n^{(k)}(1,1)$ are well known and have been extensively studied, cf. e.g. OEIS[7], A000930 and A001609. In the sequel we shall mostly set $s = 1$ and write $F_n^{(k)}(x)$ instead of $F_n^{(k)}(x,1)$.

**Remark**

As in the case of ordinary Fibonacci and Lucas polynomials there is a simple connection with matrices: Define matrices $A_k(x,s) = \left(a_k(i,j,x,s)\right)_{i,j=0}^{k-1}$ with $a_k(i,i+1,x,s) = 1$ for $0 \leq i < k-1$, $a_k(k-1,0,x,s) = s$, $a_k(k-1,k-1,x,s) = x$ and $a_k(i,j,x,s) = 0$ else.

Then for $n > -k$ the columns of $A_k^{n+k}(x,s)$ are $\left(sF_{n-j}^{(k)}(x,s), sF_{n-j+1}^{(k)}(x,s), \cdots, sF_{n-j+k-1}^{(k)}(x,s)\right)^T$ for $0 \leq j < k-1$ and $\left(F_{n+1}^{(k)}(x,s), F_{n+2}^{(k)}(x,s), \cdots, F_{n+k}^{(k)}(x,s)\right)^T$ for $j = k-1$.

This holds for $n = 1 - k$ and follows from

$$A_k(x,s)\left(F_n^{(k)}(x,s), \cdots, F_{n+k-1}^{(k)}(x,s)\right) = \left(F_{n+1}^{(k)}(x,s), \cdots, F_{n+k}^{(k)}(x,s)\right) \text{ for each } n \in \mathbb{N}.$$

These matrices satisfy $A_k^k(x,s) = xA_k^{k-1}(x,s) + sI_k$, because each entry satisfies $F_\ell^{(k)}(x,s) = xF_{\ell-1}^{(k)}(x,s) + sF_{\ell-k}^{(k)}(x,s)$ for some $\ell$.

The generalized Lucas polynomials $L_n^{(k)}(x,s)$ can be obtained as the trace

$$L_n^{(k)}(x,s) = Tr\left(A_k^n(x,s)\right).$$

For example for $k = 3$ the first terms of the sequence $A_3^n(x,s)$ are

$$\begin{pmatrix} 1 & 0 & 0 \\ 0 & 1 & 0 \\ 0 & 0 & 1 \end{pmatrix}, \begin{pmatrix} 0 & 1 & 0 \\ 0 & 0 & 1 \\ s & 0 & x \end{pmatrix}, \begin{pmatrix} 0 & 0 & 1 \\ s & 0 & x \\ sx & s & x^2 \end{pmatrix}, \begin{pmatrix} s & 0 & x \\ sx & s & x^2 \\ sx^2 & sx & s+x^3 \end{pmatrix}, \begin{pmatrix} sx & s & x^2 \\ sx^2 & sx & s+x^3 \\ s^2+sx^3 & sx^2 & 2sx+x^4 \end{pmatrix}, \cdots.$$

We have

$$\begin{pmatrix} s & 0 & x \\ sx & s & x^2 \\ sx^2 & sx & s+x^3 \end{pmatrix} = x \begin{pmatrix} 0 & 0 & 1 \\ s & 0 & x \\ sx & s & x^2 \end{pmatrix} + s \begin{pmatrix} 1 & 0 & 0 \\ 0 & 1 & 0 \\ 0 & 0 & 1 \end{pmatrix}.$$

The first terms of the sequence $F_n^{(3)}(x,s)$ are $1, x, x^2, s+x^3, 2sx+x^4, \cdots$ and those of the sequence $L_n^{(3)}(x,s)$ are $3, x, x^2, x^3+3s, x^4+4sx, \cdots$.

Another well-known representation is given by

$$F_n^{(k)}(x,s) = \det\left(f(i,j)\right)_{i,j=0}^{n-1}$$

with $f(i,i) = x$, $f(i,i+1) = 1$, $f(i,i-k+1) = (-1)^{k-1}s$ and $f(i,j) = 0$ else.

The corresponding $\ell(i,j)$ with

$$L_n^{(k)}(x,s) = \det\left(\ell(i,j)\right)_{i,j=0}^{n-1}$$

satisfy $\ell(i,j) = f(i,j)$ for $(i,j) \neq (k-1,0)$ and $\ell(k-1,0) = (-1)^{k-1}ks$.

For example

$$F_5^{(4)}(x,s) = \det\begin{pmatrix} x & 1 & 0 & 0 & 0 \\ 0 & x & 1 & 0 & 0 \\ 0 & 0 & x & 1 & 0 \\ -s & 0 & 0 & x & 1 \\ 0 & -s & 0 & 0 & x \end{pmatrix} \text{ and } L_5^{(4)}(x,s) = \det\begin{pmatrix} x & 1 & 0 & 0 & 0 \\ 0 & x & 1 & 0 & 0 \\ 0 & 0 & x & 1 & 0 \\ -4s & 0 & 0 & x & 1 \\ 0 & -s & 0 & 0 & x \end{pmatrix}.$$

## 2. Main results

**Theorem 1**

*Let*

$$A_n(k,r,x) = \left(\binom{i-k+1+r}{j}x^k + \binom{i+1+r}{j+1}\right)_{i,j=0}^{n-1} \tag{9}$$

*for $0 \leq r < k$.*

*Then*

$$x^r \det\left(A_n(k,r,x)\right) = F_{kn+r}^{(k)}(x). \tag{10}$$

**Theorem 2**

*Let*

$$B_n(k,r,x) = \left(\binom{i-k+1+r}{i-j}x^k + \binom{i+1+r}{i-j+1}\right)_{i,j=0}^{n-1} \tag{11}$$

*for $0 \leq r < k$.*

*Then*

$$x^r \det\left(B_n(k,r,x)\right) = F_{kn+r}^{(k)}(x). \tag{12}$$

**Theorem 3**

*Let*

$$C_n(k,r,x) = \left(\binom{ki+r}{i-j}x^k + \binom{k(i+1)+r}{i-j+1}\right)_{i,j=0}^{n-1} \tag{13}$$

*for $0 \leq r < k$.*

*Then*

$$x^r \det\left(C_n(k,r,x)\right) = L_{kn+r}^{(k)}(x). \tag{14}$$

**Theorem 4**

*Let*

$$D_n(k,x) = \left(\binom{i+k-1}{j+k-1}x^k - \binom{i}{j-1}\right)_{i,j=0}^{n-1}. \tag{15}$$

*Then*

$$\det(D_n(k,x)) = x^n F^{(k)}_{(k-1)n}(x). \tag{16}$$

Theorems 1 and 2 also have interesting $q-$analogues. Let us consider the generalized $q-$Fibonacci polynomials (cf. [1],[8])

$$F_n^{(k)}(x,s;q) = \sum_{j=0}^{\lfloor \frac{n}{k} \rfloor} q^{k\binom{j}{2}} \begin{bmatrix} n-(k-1)j \\ j \end{bmatrix} s^j x^{n-kj} \tag{17}$$

which satisfy the recursion $F_{n+k}^{(k)}(x,s;q) = xF_{n+k-1}^{(k)}(x,s;q) + q^n s F_n^{(k)}(x,s;q)$ with initial values $F_n^{(k)}(x,s;q) = x^n$ for $0 \leq n < k$ or equivalently with initial values $F_n^{(k)}(x,s;q) = 0$ for $-k < n < 0$ and $F_0^{(k)}(x,s;q) = 1$. Here $\begin{bmatrix} n \\ k \end{bmatrix} = \begin{bmatrix} n \\ k \end{bmatrix}_q$ denotes the $q-$binomial coefficient defined by $\begin{bmatrix} n \\ k \end{bmatrix} = \frac{(1-q^n)(1-q^{n-1})\cdots(1-q^{n-k+1})}{(1-q)(1-q^2)\cdots(1-q^k)}$ for $0 \leq k \leq n$ and $\begin{bmatrix} n \\ k \end{bmatrix} = 0$ else.

For $s=1$ we write $F_n^{(k)}(x,1;q) = F_n^{(k)}(x;q)$. For most formulas it is sufficient to choose $s=1$. But the bivariate polynomials $F_n^{(k)}(x,s;q)$ satisfy also the more elegant recursion

$$F_{n+k}^{(k)}(x,s;q) = xF_{n+k-1}^{(k)}(x,qs;q) + sF_n^{(k)}(x,q^k s;q).$$

For $k=1$ we get $F_n^{(1)}(x;q) = (x+1)(x+q)\cdots(x+q^{n-1})$.

**Theorem 5**

*Let*

$$A_n(k,r,x;q) = \left(q^{(k-r-1)j}\begin{bmatrix} i-k+r+1 \\ j \end{bmatrix}x^k + q^{(k-r)j}\begin{bmatrix} i+r+1 \\ j+1 \end{bmatrix}\right)_{i,j=0}^{n-1} \tag{18}$$

*for $0 \leq r < k$.*

*Then*

$$x^r \det\left(A_n(k,r,x;q)\right) = F^{(k)}_{kn+r}(x;q). \tag{19}$$

**Theorem 6**

*Let*

$$B_n(k,r,x;q) = \left(\begin{bmatrix} i-k+r+1 \\ i-j \end{bmatrix} x^k + q^{(k-1)j} \begin{bmatrix} i+r+1 \\ i-j+1 \end{bmatrix}\right)_{i,j=0}^{n-1} \tag{20}$$

*for $0 \leq r < k$. Then also*

$$x^r \det\left(B_n(k,r,x;q)\right) = F^{(k)}_{kn+r}(x;q). \tag{21}$$

For Theorem 3 I could only find a $q-$analogue for $k=2$:

**Proposition 7**

*Let*

$$Luc_n(x) = \sum_{k=0}^{\lfloor \frac{n}{2} \rfloor} q^{\binom{k}{2}} \begin{bmatrix} n-k \\ k \end{bmatrix} \frac{[n]}{[n-k]} x^{n-2k}. \tag{22}$$

*Then*

$$\det\left(\begin{bmatrix} 2i \\ i-j \end{bmatrix} x^2 + \begin{bmatrix} 2i+2 \\ i+1-j \end{bmatrix}\right)_{i,j=0}^{n-1} = Luc_{2n}(x),$$

$$x \det\left(\begin{bmatrix} 2i+1 \\ i-j \end{bmatrix} x^2 + \begin{bmatrix} 2i+3 \\ i+1-j \end{bmatrix}\right)_{i,j=0}^{n-1} = Luc_{2n+1}(x). \tag{23}$$

Also for Theorem 4 I have only a $q-$ analogue for $k=2$.

**Proposition 8**

Let $F_n(x,q;q) = \sum_{j=0}^{\lfloor \frac{n}{2} \rfloor} q^{j^2} \begin{bmatrix} n-j \\ j \end{bmatrix} x^{n-2j}.$

*Then*

$$\det\left(\begin{bmatrix} i+1 \\ j+1 \end{bmatrix} x^2 - \begin{bmatrix} i \\ j-1 \end{bmatrix}\right)_{i,j=0}^{n-1} = x^n F_n(x;q). \tag{24}$$

## 3. Proofs

A main tool is the following simple consequence of Cramer's rule:

**Lemma 9 (cf.[5])**

Let $T = \bigl(t(i,j)\bigr)_{i,j \geq 0}$ be a lower triangular matrix with $t(i,i) = 1$ for all $i$ and let $T_{n,1} = \bigl(t(i+1,j)\bigr)_{i,j=0}^{n-1}$. If there are numbers $M_n$ such that

$$\sum_{j=0}^{n}(-1)^{n-j}t(n,j)M_j = [n=0], \tag{25}$$

then $\det\bigl(T_{n,1}\bigr) = M_n$.

As a simple application let us derive the determinant representation of $F_n^{(k)}(x;q)$.

Let $t_k(i,i) = 1$, $t_k(i,i-1) = x$, $t_k(i,i-k) = (-1)^{k-1}q^{i-k}$ and $t_k(i,j) = 0$ else. Then (25) implies $M_0 = 1$, $M_1 = x$, $\cdots$, $M_{k-1} = x^{k-1}$, $M_k = x^k + 1$, and $M_{n+k} = xM_{n+k-1} + q^n M_n$ for $n \geq 0$. Therefore we get $M_n = F_n^{(k)}(x;q)$.

This implies $F_n^{(k)}(x;q) = \det\bigl(t_k(i+1,j)\bigr)_{i,j=0}^{n-1}$.

For example we get

$$\det\begin{pmatrix} x & 1 & 0 & 0 & 0 \\ 0 & x & 1 & 0 & 0 \\ 1 & 0 & x & 1 & 0 \\ 0 & q & 0 & x & 1 \\ 0 & 0 & q^2 & 0 & x \end{pmatrix} = F_5^{(3)}(x;q) = x^5 + \bigl(1+q+q^2\bigr)x^2.$$

If we replace $t_k(k,0)$ by $(-1)^{k-1}[k]$ we would get a $q$ – analogue of $L_n^{(k)}(x)$. But this does not have nice coefficients.

**Lemma 10**

Let $0 \leq r < k$. For $n > 0$ we get

$$\sum_{j=0}^{n}(-1)^j \binom{r+1}{j} F_{k(n-j)+r}^{(k)}(x) = x^{r+1} F_{kn-1}^{(k)}(x) \tag{26}$$

and for $0 < i < k+n$

$$\sum_{j=0}^{n-1}(-1)^{n-1-j}\binom{i-k+r}{n-1-j}x^k F_{kj+r}^{(k)}(x) = x^{r+i}F_{kn-i}^{(k)}(x). \qquad (27)$$

**Remark**

If we consider the shift operator $E$ defined by $EF_n^{(k)} = F_{n+1}^{(k)}$ then the recursion of $F_n^{(k)}$ can symbolically be written as $\left(E^k - 1\right)F_n^{(k)}(x) = xE^{k-1}F_n^{(k)}(x)$ or $\left(1 - E^{-k}\right)F_n^{(k)}(x) = xE^{-1}F_n^{(k)}(x)$ and the identities (26) as

$$\left(1 - E^{-k}\right)^{r+1} F_{kn+r}(x) = \left(xE^{-1}\right)^{r+1} F_{kn+r}(x) = x^{r+1}F_{kn-1}(x)$$

and (27) as

$$\left(1 - E^{-k}\right)^{i+r-k} F_{k(n-1)+r}^{(k)}(x) = \left(xE^{-1}\right)^{i+r-k} F_{k(n-1)+r}^{(k)}(x) = x^{i+r-k}F_{k(n-1)+r-i-r+k}^{(k)} = x^{i+r-k}F_{kn-i}^{(k)}(x).$$

I shall not try to make this exact but give a proof by induction.

1) Identity (26).

For $r = 0$ and $n > 0$ it reduces to the defining recursion $F_{kn}^{(k)}(x) - F_{k(n-1)}^{(k)}(x) = xF_{kn-1}^{(k)}(x)$.

Suppose it holds for $r - 1$ and $n > 0$. Then we get for $r \leq n$

$$x^{r+1}F_{kn-1}^{(k)}(x) = xx^r F_{kn-1}^{(k)}(x) = \sum_{j=0}^{n}(-1)^j\binom{r}{j}xF_{k(n-j)+r-1}^{(k)}(x)$$

$$= \sum_{j=0}^{n}(-1)^j\binom{r}{j}\left(F_{k(n-j)+r}^{(k)}(x) - F_{k(n-j-1)+r}^{(k)}(x)\right) = \sum_{j=0}^{n}(-1)^j\binom{r}{j}F_{k(n-j)+r}^{(k)}(x) - \sum_{j=0}^{n}(-1)^j\binom{r}{j}F_{k(n-j-1)+r}^{(k)}(x)$$

$$= \sum_{j=0}^{n}(-1)^j\binom{r}{j}F_{k(n-j)+r}^{(k)}(x) + \sum_{j=1}^{n+1}(-1)^j\binom{r}{j-1}F_{k(n-j)+r}^{(k)}(x)$$

$$= \sum_{j=0}^{n}(-1)^j\binom{r+1}{j}F_{k(n-j)+r}^{(k)}(x) + (-1)^{n+1}\binom{r}{n}F_{-k+r}^{(k)}(x).$$

The last term vanishes because $F_{-k+r}^{(k)}(x) = 0$.

Let now $n \leq r < k$. Since the identity is true for $r = n$ we get for $r + 1$

$$x^{r+1}F_{kn-1}^{(k)}(x) = xx^r F_{kn-1}^{(k)}(x) = \sum_{j=0}^{n-1}(-1)^j\binom{r}{j}xF_{k(n-j)+r-1}^{(k)}(x) + (-1)^n\binom{r}{n}xF_{r-1}^{(k)}(x)$$

$$= \sum_{j=0}^{n-1}(-1)^j\binom{r}{j}\left(F_{k(n-j)+r}^{(k)}(x) - F_{k(n-j-1)+r}^{(k)}(x)\right) + (-1)^n\binom{r}{n}xF_{r-1}^{(k)}(x)$$

$$= \sum_{j=0}^{n-1}(-1)^j\binom{r}{j}F_{k(n-j)+r}^{(k)}(x) + \sum_{j=1}^{n}(-1)^j\binom{r}{j-1}F_{k(n-j)+r}^{(k)}(x) + (-1)^n\binom{r}{n}xF_{r-1}^{(k)}(x)$$

$$= \sum_{j=0}^{n-1}(-1)^j \binom{r+1}{j} F_{k(n-j)+r}^{(k)}(x) + (-1)^n \binom{r}{n-1} F_r^{(k)}(x) + (-1)^n x \binom{r}{n} F_{r-1}^{(k)}(x)$$

$$= \sum_{j=0}^{n-1}(-1)^j \binom{r+1}{j} F_{k(n-j)+r}^{(k)}(x) + (-1)^n \binom{r}{n-1} x^r + (-1)^n \binom{r}{n} x^r$$

$$= \sum_{j=0}^{n}(-1)^j \binom{r+1}{j} F_{k(n-j)+r}^{(k)}(x).$$

For $n = 0$ the left-hand side of (26) is $x^r$.

2) Identity (27).

Let $s(n,i) = \sum_{j=0}^{n-1}(-1)^{n-1-j} \binom{i-k+r}{n-1-j} x^k F_{kj+r}^{(k)}(x)$ and $f(n,i) = x^{r+i} F_{kn-i}^{(k)}(x)$.

We show first that for $0 \leq r < k$ and $n > 0$ the identity $s(n,i) = f(n,i)$ holds for $0 < i \leq k - r$.

For example for $(k,r,i) = (2,0,1)$ we get $xF_{2n-1}(x) = x^2 \sum_{j=0}^{n-1} F_{2j}(x)$.

We have $s(n, k-r) = f(n, k-r)$ because $x^{r+i} F_{kn-i}^{(k)}(x) = x^k F_{k(n-1)+r}^{(k)}$.

Suppose we know the formula for $i + 1 \leq k - r$. We want to prove it for $i$ and all $n$.

For $n = 1$ formula (27) reduces to $s(1,i) = x^{r+i} F_{k-i}^{(k)}(x) = x^{r+k} = x^k F_r^{(k)}(x) = f(1,i)$

Suppose it holds for $n - 1$.

Then we get from $x^{r+i} F_{kn-i}^{(k)}(x) = x^{r+i+1} F_{kn-i-1}^{(k)}(x) + x^{r+i} F_{k(n-1)-i}^{(k)}(x)$

$$f(n,i) = x^{r+i} F_{kn-i}^{(k)}(x) = \sum_{j=0}^{n-1}(-1)^{n-1-j} \binom{i+1-k+r}{n-1-j} x^k F_{kj+r}^{(k)}(x) - \sum_{j=0}^{n-1}(-1)^{n-1-j} \binom{i-k+r}{n-2-j} x^k F_{kj+r}^{(k)}(x)$$

$$= \sum_{j=0}^{n-1}(-1)^{n-1-j} \binom{i-k+r}{n-1-j} x^k F_{kj+r}^{(k)}(x) = s(n,i).$$

For $i = 1$ this gives

$$x^{r+1} F_{kn-1}^{(k)}(x) = \sum_{j=0}^{n-1}(-1)^{n-1-j} \binom{1-k+r}{n-1-j} x^k F_{kj+r}^{(k)}(x). \tag{28}$$

In the other direction we have for $i > k - r$

$$s(n, i+1) = s(n,i) - s(n-1, i) \tag{29}$$

because

$$s(n,i) - s(n-1,i) = \sum_{j=0}^{n-1}(-1)^{n-1-j}\binom{i-k+r}{n-1-j}x^k F_{kj+r}^{(k)}(x) + \sum_{j=0}^{n-2}(-1)^{n-1-j}\binom{i-k+r}{n-2-j}x^k F_{kj+r}^{(k)}(x)$$

$$= \sum_{j=0}^{n-2}(-1)^{n-1-j}\left[\binom{i-k+r}{n-1-j}+\binom{i-k+r}{n-2-j}\right]x^k F_{kj+r}^{(k)}(x) + \binom{i-k+r}{0}x^k F_{k(n-1)+r}^{(k)}(x)$$

$$= \sum_{j=0}^{n-1}(-1)^{n-1-j}\binom{i+1-k+r}{n-1-j}x^k F_{kj+r}^{(k)}(x) = s(n,i+1).$$

On the other hand the recursion $x^{r+i}F_{kn-i}^{(k)}(x) = x^{r+i+1}F_{kn-i-1}^{(k)}(x) + x^{r+i}F_{k(n-1)-i}^{(k)}(x)$ gives

$$f(n,i+1) = f(n,i) - f(n-1,i), \qquad (30)$$

but only for $i < kn$, because

$$f(n,kn) = x^{r+kn}F_{kn-kn}^{(k)}(x) = x^{r+kn} \neq f(n,kn+1) + f(n-1,kn) = 0.$$

Comparing (29) and (30) we get

$$s(n,i+1) - f(n,i+1) = s(n,i) - f(n,i) - \big(s(n-1,i) - f(n-1,i)\big). \qquad (31)$$

Since $s(1,k+1) = x^{k+r} \neq f(1,k+1) = 0$ we see by induction that $s(n,i) = f(n,i)$ for all $i < n+k$, but that $s(n,n+k) \neq f(n,n+k)$.

**Proof of Theorem 1**

For $k=1$ the matrix $A_n(1,0,x)$ is triangular with $1+x$ in the main diagonal. Therefore $\det A_n(1,0,x) = (1+x)^n$.

For arbitrary $k$ we get

$$\det\left(\binom{i-k+1+r}{j}x^k + \binom{i+1+r}{j+1}\right)_{i,j=0}^{n-1} = \det\left(\binom{r+1-k}{i-j}x^k + \binom{r+1}{i+1-j}\right)_{i,j=0}^{n-1} \qquad (32)$$

because

$$\left(\binom{i}{j}\right)_{i,j=0}^{n-1} \cdot \left(\binom{r+1}{j-i+1}\right)_{i,j=0}^{n-1} = \left(\binom{i+1+r}{j+1}\right)_{i,j=0}^{n-1} \quad \text{and}$$

$$\left(\binom{i}{j}\right)_{i,j=0}^{n-1} \cdot \left(\binom{r+1-k}{j-i}\right)_{i,j=0}^{n-1} = \left(\binom{i-k+1+r}{j}\right)_{i,j=0}^{n-1}.$$

This follows from Vandermonde's formula

$$\sum_j \binom{i}{j}\binom{r+1}{\ell-j+1} = \binom{i+r+1}{\ell+1}.$$

Combining these results and observing that the inverse of $\left(\binom{i}{j}\right)_{i,j=0}^{n-1}$ is $\left((-1)^{i-j}\binom{i}{j}\right)_{i,j=0}^{n-1}$

we get

$$\left((-1)^{i-j}\binom{i}{j}\right)_{i,j=0}^{n-1} A_n(k,r,x) = \left(\binom{r+1-k}{j-i}x^k + \binom{r+1}{j-i+1}\right)_{i,j=0}^{n-1}. \tag{33}$$

Transposing this matrix we get (32).

In the matrices

$$\left(\binom{r+1-k}{i-j}x^k + \binom{r+1}{i-j+1}\right)_{i,j=0}^{n-1}$$

all entries for $j > i+1$ vanish. For $j = i+1$ all entries are 1. By Lemma 9 this is equivalent with the identity

$$\sum_{j=0}^{n}(-1)^{n-j}t(n,j)\frac{F_{kj+r}^{(k)}(x)}{x^r} = [n=0]. \tag{34}$$

for $t(i,j) = \binom{r+1}{i-j} + x^k\binom{r+1-k}{i-j-1}$.

Combining (26) and (28) we get (34). Note that for $n=0$ the left-hand side of (26) is $x^r$ and the left-hand side of (28) is 0.

**Proof of Theorem 2**

The entries $b(i,j)$ of the matrix $B_n(k,r,x)$ satisfy $b(i,j) = 0$ for $j > i+1$ and $b(i,i+1) = 1$. Thus we can apply Lemma 9 with

$$t(i,j) = \binom{i+r-k}{i-j-1}x^k + \binom{i+r}{i-j}.$$

We must show that

$$\sum_{j=0}^{i} t(i,j)(-1)^{i-j}F_{kj+r}^{(k)}(x) = x^r[i=0].$$

It suffices to show that for $0 \leq r < k$

$$\sum_{j=0}^{n}(-1)^{n-j}\binom{n+r}{n-j}F_{kj+r}^{(k)}(x) = x^{n+r}F_{(k-1)n}^{(k)}(x) \tag{35}$$

and

$$\sum_{j=0}^{n}(-1)^{n-j}\binom{n+r-k}{n-j-1}x^k F_{kj+r}^{(k)}(x) = -x^{n+r}F_{(k-1)n}^{(k)}(x). \tag{36}$$

Both formulae are a special case of (27).

**Proof of Theorem 3**

$$C_n(k,r,x) = \left(\binom{ki+r}{i-j}x^k + \binom{k(i+1)+r}{i-j+1}\right)_{i,j=0}^{n-1} \text{ is a Hessenberg matrix with 1 in the diagonal}$$

$j = i+1$. We want to show that

$$\det\left(\binom{ki+r}{i-j}x^k + \binom{k(i+1)+r}{i-j+1}\right)_{i,j=0}^{n-1} = \frac{L_{kn+r}^{(k)}}{x^r}. \tag{37}$$

Let now $l_n^{(k)}(x) = L_n^{(k)}(x)$ for $n > 0$ and $l_0^{(k)}(x) = 1$.

(37) holds for $0 \leq r < k$ if

$$\sum_{j=0}^{n}(-1)^{n-j}\left(\binom{k(n-1)+r}{n-j-1}x^k + \binom{kn+r}{n-j}\right)\frac{l_{kj+r}^{(k)}}{x^r} = [n=0]. \tag{38}$$

We know that (cf. [3]) $\sum_{j=0}^{\lfloor\frac{n}{k}\rfloor}(-1)^j\binom{n}{j}l_{n-kj}^{(k)}(x) = x^n$.

Thus for $n > 0$ and $r < k$ we have

$$\sum_{j=0}^{n}(-1)^{n-j}\binom{kn+r}{n-j}l_{kj+r}^{(k)}(x) = x^{kn+r} \text{ and}$$

$$\sum_{j=0}^{n}(-1)^{n-1-j}\binom{k(n-1)+r}{n-1-j}l_{kj+r}^{(k)}(x) = x^{k(n-1)+r}$$

For $n = 0$ the first sum is $x^r$ and the second sum vanishes. Thus (38) is true.

**Proof of Theorem 4**

By Lemma 9 it suffices to show that

$$\sum_{j=0}^{n-1} \binom{n+k-2}{n-j-1} x^{k+j} F^{(k)}_{(k-1)j}(x) = \sum_{j=0}^{n} \binom{n-1}{n-j} x^j F^{(k)}_{(k-1)j}(x)$$

for $n > 0$. We show that these sums are equal to $xF^{(k)}_{kn-1}(x)$. This follows from

$$\sum_{j=0}^{n-1} \binom{n+k-2}{n-j-1} x^{k+j} F^{(k)}_{(k-1)j}(x) = \sum_j \binom{n+k-2}{n-j-1} x^{k+j} \sum_\ell \binom{(k-1)j-(k-1)\ell}{\ell} x^{(k-1)j-k\ell}$$

$$= \sum_r x^{kn-kr} \sum_{j-\ell=n-r-1} \binom{n+k-2}{n-j-1}\binom{(k-1)j-(k-1)\ell}{\ell}$$

$$= \sum_r x^{kn-kr} \sum_\ell \binom{n+k-2}{r-\ell}\binom{(k-1)(n-r-1)}{\ell} = \sum_r x^{kn-kr} \binom{n+k-2+kn-n-kr+r-k+1}{r}$$

$$= \sum_r x^{kn-kr} \binom{kn-kr+r-1}{r} = xF^{(k)}_{kn-1}(x).$$

$$\sum_{j=0}^{n-1} \binom{n-1}{n-j} x^j F^{(k)}_{(k-1)j}(x) = \sum_j \binom{n-1}{n-j} x^j \sum_\ell \binom{(k-1)j-(k-1)\ell}{\ell} x^{(k-1)j-k\ell}$$

$$= \sum_r x^{kn-kr} \sum_{j-\ell=n-r} \binom{n-1}{n-j-1}\binom{(k-1)j-(k-1)\ell}{\ell}$$

$$= \sum_r x^{kn-kr} \sum_\ell \binom{n-1}{r-\ell}\binom{(k-1)j-(k-1)\ell}{\ell} = \sum_r x^{kn-kr} \binom{kn-1-kr+r}{r} = xF^{(k)}_{kn-1}(x).$$

As $q-$analogue of Lemma 10 we get

**Lemma 11**

Let $0 \leq r < k$. For $n > 0$ we get

$$\sum_{j=0}^{n}(-1)^j \begin{bmatrix} r+1 \\ j \end{bmatrix} q^{nkj+\binom{j}{2}-k\binom{j+1}{2}} F^{(k)}_{k(n-j)+r}(x;q) = x^{r+1} F^{(k)}_{kn-1}(x;q) \tag{39}$$

and for $0 < i < k+n$

$$\sum_{j=0}^{n-1}(-1)^{n-1-j} \begin{bmatrix} i-k+r \\ n-1-j \end{bmatrix} x^k q^{\binom{n-j-1}{2}-k\binom{n-j}{2}+(n-j-1)(kn-i+1)} F^{(k)}_{kj+r}(x;q) = x^{r+i} F^{(k)}_{kn-i}(x;q). \tag{40}$$

**Proof**

Let us first prove (39). For $r = 0$ it reduces to the defining identity

$$\sum_{j=0}^{n}(-1)^{j}\begin{bmatrix}1\\j\end{bmatrix}q^{nkj+\binom{j}{2}-k\binom{j+1}{2}}F_{k(n-j)}^{(k)}(x;q)=F_{kn}^{(k)}(x;q)-q^{k(n-1)}F_{k(n-1)}^{(k)}(x;q)=xF_{kn-1}^{(k)}(x;q).$$

Suppose that (39) holds for $r-1$ and $n>0$. Then we get for $r \leq n$

$$x^{r+1}F_{kn-1}^{(k)}(x;q)=xx^{r}F_{kn-1}^{(k)}(x;q)=\sum_{j=0}^{n}(-1)^{j}\begin{bmatrix}r\\j\end{bmatrix}q^{nkj+\binom{j}{2}-k\binom{j+1}{2}}xF_{k(n-j)+r-1}^{(k)}(x;q)$$

$$=\sum_{j=0}^{n}(-1)^{j}\begin{bmatrix}r\\j\end{bmatrix}q^{nkj+\binom{j}{2}-k\binom{j+1}{2}}\left(F_{k(n-j)+r}^{(k)}(x;q)-q^{k(n-j-1)+r}F_{k(n-j-1)+r}^{(k)}(x;q)\right)$$

$$=\sum_{j=0}^{n}(-1)^{j}\begin{bmatrix}r\\j\end{bmatrix}q^{nkj+\binom{j}{2}-k\binom{j+1}{2}}F_{k(n-j)+r}^{(k)}(x;q)+\sum_{j=1}^{n+1}(-1)^{j}\begin{bmatrix}r\\j-1\end{bmatrix}q^{nk(j-1)+\binom{j-1}{2}-k\binom{j}{2}}q^{k(n-j)+r}F_{k(n-j)+r}^{(k)}(x;q)$$

$$=\sum_{j=0}^{n}(-1)^{j}\left(\begin{bmatrix}r\\j\end{bmatrix}+q^{r-j+1}\begin{bmatrix}r\\j-1\end{bmatrix}\right)q^{nkj+\binom{j}{2}-k\binom{j+1}{2}}F_{k(n-j)+r}^{(k)}(x;q)+(-1)^{n+1}\begin{bmatrix}r\\n\end{bmatrix}q^{nk(n+1)+\binom{n+1}{2}-k\binom{n+2}{2}}F_{-k+r}^{(k)}(x;q)$$

$$=\sum_{j=0}^{n}(-1)^{j}\begin{bmatrix}r+1\\j\end{bmatrix}q^{nkj+\binom{j}{2}-k\binom{j+1}{2}}F_{k(n-j)+r}^{(k)}(x;q)$$

Note that for $n>r$ the coefficient $\begin{bmatrix}r\\n\end{bmatrix}$ vanishes and for $r=n$ we have $F_{-k+r}^{(k)}(x)=0$ because $r<k$.

There remains to study the case $n \leq r < k$. Since (49) is true for $r=n$ we get by induction

$$x^{r+1}F_{kn-1}^{(k)}(x;q)=xx^{r}F_{kn-1}^{(k)}(x;q)=\sum_{j=0}^{n-1}(-1)^{j}\begin{bmatrix}r\\j\end{bmatrix}q^{nkj+\binom{j}{2}-k\binom{j+1}{2}}xF_{k(n-j)+r-1}^{(k)}(x;q)+(-1)^{n}\begin{bmatrix}r\\n\end{bmatrix}q^{n^{2}k+\binom{n}{2}-k\binom{n+1}{2}}xF_{r-1}^{(k)}(x;q)$$

$$=\sum_{j=0}^{n-1}(-1)^{j}\begin{bmatrix}r\\j\end{bmatrix}q^{nkj+\binom{j}{2}-k\binom{j+1}{2}}\left(F_{k(n-j)+r}^{(k)}(x;q)-q^{k(n-j-1)+r}F_{k(n-j-1)+r}^{(k)}(x;q)\right)+(-1)^{n}\begin{bmatrix}r\\n\end{bmatrix}q^{n^{2}k+\binom{n}{2}-k\binom{n+1}{2}}xF_{r-1}^{(k)}(x;q)$$

$$=\sum_{j=0}^{n-1}(-1)^{j}\begin{bmatrix}r\\j\end{bmatrix}q^{nkj+\binom{j}{2}-k\binom{j+1}{2}}F_{k(n-j)+r}^{(k)}(x;q)+\sum_{j=1}^{n-1}(-1)^{j}\begin{bmatrix}r\\j-1\end{bmatrix}q^{nk(j-1)+\binom{j-1}{2}-k\binom{j}{2}+k(n-j)+r}F_{k(n-j)+r}^{(k)}(x;q)$$

$$+(-1)^{n}\begin{bmatrix}r\\n-1\end{bmatrix}q^{nk(n-1)+\binom{n-1}{2}-k\binom{n}{2}}q^{r}F_{r}^{(k)}(x;q)+(-1)^{n}\begin{bmatrix}r\\n\end{bmatrix}q^{n^{2}k+\binom{n}{2}-k\binom{n+1}{2}}xF_{r-1}^{(k)}(x;q)$$

$$=\sum_{j=0}^{n}(-1)^{j}\begin{bmatrix}r+1\\j\end{bmatrix}q^{nkj+\binom{j}{2}-k\binom{j+1}{2}}F_{k(n-j)+r}^{(k)}(x;q)$$

because

$$F_{r}^{(k)}(x;q)=xF_{r-1}^{(k)}(x;q).$$

Now we show (40) for $1 \leq i \leq k - r$.

Let $s(n, i, q) = \sum_{j=0}^{n-1}(-1)^{n-1-j} \begin{bmatrix} i-k+r \\ n-1-j \end{bmatrix} x^k q^{\binom{n-j-1}{2}-k\binom{n-j}{2}+(n-j-1)(kn-i+1)} F_{kj+r}^{(k)}(x;q)$ and

$f(n, i, q) = x^{r+i} F_{kn-i}^{(k)}(x; q).$

We show first that for $0 \leq r < k$ and $n > 0$ the identity $s(n, i) = f(n, i)$ holds for $0 < i \leq k - r$.

We have $s(n, k - r, q) = f(n, k - r, q)$ because
$f(n, k-r, q) = x^{r+k-r} F_{kn+r-k}^{(k)}(x;q) = x^k F_{k(n-1)+r}^{(k)}(x;q) = s(n, k-r, q).$

Suppose we know the formula for $i + 1 \leq k - r$. We want to prove it for $i$ and all $n$.

For $n = 1$ formula (27) reduces to
$s(1, i, q) = x^{r+i} F_{k-i}^{(k)}(x; q) = x^{r+i} x^{k-i} = x^{k+r} = x^k F_r^{(k)}(x; q) = f(1, i, q).$

Suppose it holds for $n - 1$.

Then we get from $x^{r+i} F_{kn-i}^{(k)}(x;q) = x^{r+i+1} F_{kn-i-1}^{(k)}(x;q) + q^{k(n-1)-i} x^{r+i} F_{k(n-1)-i}^{(k)}(x;q)$

$f(n, i, q) = x^{r+i} F_{kn-i}^{(k)}(x;q) = \sum_{j=0}^{n-1}(-1)^{n-1-j} \begin{bmatrix} r-k+i+1 \\ n-1-j \end{bmatrix} x^k F_{kj+r}^{(k)}(x;q) q^{j+1-n} q^{\binom{n-j-1}{2}-k\binom{n-j}{2}+(n-j-1)(kn-i)}$

$- q^{k(n-1)-i} \sum_{j=0}^{n-2}(-1)^{n-1-j} \begin{bmatrix} r-k+i \\ n-2-j \end{bmatrix} x^k F_{kj+r}^{(k)}(x;q) q^{\binom{n-j-2}{2}-k\binom{n-j-1}{2}+(n-j-2)(kn-i+1-k)}$

$= \sum_{j=0}^{n-1}(-1)^{n-1-j} x^k F_{kj+r}^{(k)}(x;q) q^{\binom{n-j-1}{2}-k\binom{n-j}{2}+(n-j-1)(kn-i+1)} q^{j+1-n} \left( \begin{bmatrix} r-k+i+1 \\ n-1-j \end{bmatrix} - \begin{bmatrix} r-k+i \\ n-2-j \end{bmatrix} \right)$

$= \sum_{j=0}^{n-1}(-1)^{n-1-j} x^k F_{kj+r}^{(k)}(x;q) q^{\binom{n-j-1}{2}-k\binom{n-j}{2}+(n-j-1)(kn-i+1)} \begin{bmatrix} r-k+i \\ n-1-j \end{bmatrix} = s(n, i, q).$

To show the other direction we first prove

$$s(n, i+1, q) = s(n, i, q) - q^{k(n-1)-i} s(n-1, i, q). \tag{41}$$

This follows from

$s(n, i, q) - q^{k(n-1)-i} s(n-1, i, q) = \sum_{j=0}^{n-1}(-1)^{n-1-j} \begin{bmatrix} i-k+r \\ n-1-j \end{bmatrix} x^k q^{\binom{n-j-1}{2}-k\binom{n-j}{2}+(n-j-1)(kn-i+1)} F_{kj+r}^{(k)}(x;q)$

$- q^{k(n-1)-i} \sum_{j=0}^{n-1}(-1)^{n-2-j} \begin{bmatrix} i-k+r \\ n-2-j \end{bmatrix} x^k q^{\binom{n-j-2}{2}-k\binom{n-j-1}{2}+(n-j-2)(kn-i+1-k)} F_{kj+r}^{(k)}(x;q)$

$= \sum_{j=0}^{n-1}(-1)^{n-1-j} q^{n-j-1} \begin{bmatrix} i-k+r \\ n-1-j \end{bmatrix} x^k q^{\binom{n-j-1}{2}-k\binom{n-j}{2}+(n-j-1)(kn-i)} F_{kj+r}^{(k)}(x;q)$

$$+\sum_{j=0}^{n-1}(-1)^{n-1-j}\begin{bmatrix}i-k+r\\n-2-j\end{bmatrix}x^k q^{\binom{n-j-1}{2}-k\binom{n-j}{2}+(n-j-1)(kn-i)}F_{kj+r}^{(k)}(x;q)$$

$$=\sum_{j=0}^{n-1}(-1)^{n-1-j}\begin{bmatrix}i+1-k+r\\n-1-j\end{bmatrix}x^k q^{\binom{n-j-1}{2}-k\binom{n-j}{2}+(n-j-1)(kn-i)}F_{kj+r}^{(k)}(x;q)=s(n,i+1,q).$$

On the other hand the recursion $x^{r+i}F_{kn-i}^{(k)}(x;q) = x^{r+i+1}F_{kn-i-1}^{(k)}(x;q) + q^{k(n-1)-i}x^{r+i}F_{k(n-1)-i}^{(k)}(x;q)$ gives

$$f(n, i+1, q) = f(n, i, q) - q^{k(n-1)-i} f(n-1, i, q) \tag{42}$$

for $i < kn$.

Comparing (41) and (42) we get

$$s(n,i+1,q) - f(n,i+1,q) = s(n,i,q) - f(n,i,q) - q^{k(n-1)-i}\Big(s(n-1,i,q) - f(n-1,i,q)\Big).$$

Since $s(1, k+1, q) = x^{k+r} \neq f(1, k+1, q) = 0$ we see by induction that
$s(n, n+k, q) \neq f(n, n+k, q)$.

**Proof of Theorem 4**

Recall that

$$A_n(k,r,x;q) = \left(q^{(k-r-1)j}\begin{bmatrix}i-k+r+1\\j\end{bmatrix}x^k + q^{(k-r)j}\begin{bmatrix}i+r+1\\j+1\end{bmatrix}\right)_{i,j=0}^{n-1}$$

We show first that

$$\left((-1)^{i-j}q^{\binom{i-j}{2}}\begin{bmatrix}i\\j\end{bmatrix}\right)_{i,j=0}^{n-1}A_n(k,r,x;q) = \Big(h(i,j,k,r)\Big)_{i,j=0}^{n-1} \tag{43}$$

where $h(i, j, k, r) = 0$ for $j < i - 1$ and

$$h(i,j,k,r) = q^{jk+(i-j)(i+r)}\begin{bmatrix}r+1\\j-i+1\end{bmatrix} + q^{(j-i)(k-r-i-1)}\begin{bmatrix}r-k+1\\j-i\end{bmatrix}x^k$$

else.

The identity (43) is equivalent with

$$A_n(k,r,x;q) = \left(q^{\binom{i-j}{2}}\begin{bmatrix}i\\j\end{bmatrix}\right)_{i,j=0}^{n-1}\Big(h(i,j,k,r)\Big)_{i,j=0}^{n-1}. \tag{44}$$

This is again equivalent with

$$\sum_{\ell} \begin{bmatrix} i \\ \ell \end{bmatrix} q^{(j-\ell)(k-r-\ell-1)} \begin{bmatrix} r-k+1 \\ j-\ell \end{bmatrix} = q^{j(k-r-1)} \begin{bmatrix} r+1-k+i \\ j \end{bmatrix} \quad (45)$$

and

$$\sum_{\ell} \begin{bmatrix} i \\ \ell \end{bmatrix} \begin{bmatrix} r+1 \\ j-\ell+1 \end{bmatrix} q^{(jk)+(\ell-j)(\ell+r)} = q^{(k-r)j} \begin{bmatrix} r+i+1 \\ j+1 \end{bmatrix}. \quad (46)$$

Since $(j-\ell)(k-r-\ell-1) = \ell\bigl(r-k+1-j+\ell\bigr) + j\bigl(k-r-1\bigr)$ (45) follows from $q-$Vandermonde's formula

$$\begin{bmatrix} n+m \\ k \end{bmatrix} = \sum_{j} \begin{bmatrix} n \\ j \end{bmatrix} \begin{bmatrix} n \\ k-j \end{bmatrix} q^{(n-j)(k-j)}.$$

In the same way we get

$$\sum_{\ell} \begin{bmatrix} i \\ \ell \end{bmatrix} \begin{bmatrix} r+1 \\ j-\ell+1 \end{bmatrix} q^{(jk)+(\ell-j)(\ell+r)} = q^{(k-r)j} \begin{bmatrix} r+i+1 \\ j+1 \end{bmatrix}$$

because $(jk)+(\ell-j)(\ell+r)-j(k-r) = \ell(r+\ell-j).$

In order to apply Lemma 9 we consider the transpose $\bigl(h(j,i,k,r)\bigr)_{i,j=0}^{n-1}$ and divide each row $i$ by its entry $q^{ik+i+r+1}$ for $j=i+1$ and then change $i \to i-1$ and obtain

$t(k,r,0,j) = [j=0],$

$$t(k,r,i,j) = q^{(j-i+1)(j+r)-(i+r)} \begin{bmatrix} r+1 \\ i-j \end{bmatrix} + x^{k} \begin{bmatrix} r-k+1 \\ i-j-1 \end{bmatrix} q^{(i-j-1)(k-r-j-1)-k(i-1)-i-r}.$$

Then Theorem 4 is equivalent with

$$x^{r} \det\bigl(t(k,r,i+1,j)\bigr)_{i,j=0}^{n-1} = \frac{F_{kn+r}^{(k)}(x)}{q^{nr+\binom{n+1}{2}+k\binom{n}{2}}}. \quad (47)$$

Therefore we must show that

$$\sum_{j=0}^{i}(-1)^{i-j} t(k,r,i,j) \frac{F_{kj+r}^{(k)}(x)}{x^{r} q^{jr+\binom{j+1}{2}+k\binom{j}{2}}} = [i=0]. \quad (48)$$

The coefficient of $\begin{bmatrix} r+1 \\ i-j \end{bmatrix}$ is $q^{(j-i+1)(j+r)-(i+r)-\left(jr+\binom{j+1}{2}+k\binom{j}{2}\right)}$. If we change $j \to i-j$ we get

$$\frac{q^{-k\binom{i}{2}-\binom{i+1}{2}-ir}}{x^r}\sum_{j=0}^{i}(-1)^j q^{\binom{j}{2}-k\binom{j+1}{2}+ijk}\begin{bmatrix}r+1\\j\end{bmatrix}F^{(k)}_{k(i-j)+r}(x).$$

If we do the same with the coefficient of $\begin{bmatrix}r-k+1\\i-j-1\end{bmatrix}$ we get

$$x^{k-r}\sum_{j=0}^{i}(-1)^j q^{(j-1)(k-r-i+j-1)-k(i-1)-i-r-\left[(i-j)r+\binom{i-j+1}{2}+k\binom{i-j}{2}\right]}\begin{bmatrix}r-k+1\\j-1\end{bmatrix}F^{(k)}_{k(i-j)+r}(x)$$

$$=x^{k-r}q^{-\binom{i+1}{2}-\binom{i}{2}k-ir}\sum_{j=0}^{i}(-1)^j q^{\binom{j-1}{2}-k\binom{j}{2}+ik(j-1)}\begin{bmatrix}r-k+1\\j-1\end{bmatrix}F^{(k)}_{k(i-j)+r}(x)$$

Thus it suffices to show that for $n > 0$ and $0 \leq r < k$

$$\sum_{j=0}^{n}(-1)^j \begin{bmatrix}r+1\\j\end{bmatrix} q^{nkj+\binom{j}{2}-k\binom{j+1}{2}} F^{(k)}_{k(n-j)+r}(x) = x^{r+1}F^{(k)}_{kn-1}(x) \tag{49}$$

and

$$\sum_{j=0}^{n}(-1)^j x^k q^{\binom{j-1}{2}+nk(j-1)-k\binom{j}{2}}\begin{bmatrix}r-k+1\\j-1\end{bmatrix}F^{(k)}_{k(n-j)+r}(x) = -x^{r+1}F^{(k)}_{nk-1}(x). \tag{50}$$

This follows from Lemma 11.

**Proof of Theorem 5**

$$B_n(k,r,x;q) = \left(\begin{bmatrix}i-k+r+1\\i-j\end{bmatrix}x^k + q^{(k-1)j}\begin{bmatrix}i+r+1\\i-j+1\end{bmatrix}\right)_{i,j=0}^{n-1}$$

$$x^r \det B_n(k,r,x;q) = F^{(k)}_{kn+r}(x;q).$$

In this case we have $q^{(k-1)j}$ in the diagonal $j = i+1$. If we divide row $i$ by $q^{(k-1)(i+1)}$ the new determinant is $\dfrac{F^{(k)}_{kn+r}(x)}{q^{(k-1)\binom{n+1}{2}}}$.

Therefore it suffices to show that

$$\sum_{j=0}^{n}(-1)^{n-j}\begin{bmatrix}n+r\\n-j\end{bmatrix}q^{(k-1)\left(\binom{n}{2}-\binom{j}{2}\right)}F^{(k)}_{kj+r}(x;q)=x^{n+r}F^{(k)}_{(k-1)n}(x;q) \tag{51}$$

and

$$\sum_{j=0}^{n}(-1)^{n-j}x^k\begin{bmatrix}n+r-k\\n-j-1\end{bmatrix}q^{(k-1)\left(\binom{n}{2}-\binom{j+1}{2}\right)}F^{(k)}_{kj+r}(x;q)=x^{n+r}F^{(k)}_{(k-1)n}(x;q). \tag{52}$$

Formula (52) is (40) for $i=n$.

To prove (51) let

$$h(n,k,r)=\sum_{j=0}^{n}(-1)^{n-j}\begin{bmatrix}n+r\\n-j\end{bmatrix}q^{(k-1)\left(\binom{n}{2}-\binom{j}{2}\right)}F^{(k)}_{kj+r}(x;q).$$

We show first that $xh(n,k,r)=h(n,k,r+1)$ if $r+1<k$.

$$xh(n,k,r)=\sum_{j=0}^{n}(-1)^{n-j}\begin{bmatrix}n+r\\n-j\end{bmatrix}q^{(k-1)\left(\binom{n}{2}-\binom{j}{2}\right)}xF^{(k)}_{kj+r}(x;q)$$

$$=(-1)^n\begin{bmatrix}n+r\\n\end{bmatrix}q^{(k-1)\binom{n}{2}}xF^{(k)}_r(x;q)+\sum_{j=1}^{n}(-1)^{n-j}\begin{bmatrix}n+r\\n-j\end{bmatrix}q^{(k-1)\left(\binom{n}{2}-\binom{j}{2}\right)}\left(F^{(k)}_{kj+r+1}(x;q)-q^{k(j-1)+r+1}F^{(k)}_{k(j-1)+r+1}(x;q)\right)$$

$$=(-1)^n\begin{bmatrix}n+r\\n\end{bmatrix}q^{(k-1)\binom{n}{2}}xF^{(k)}_r(x;q)+\sum_{j=1}^{n}(-1)^{n-j}\begin{bmatrix}n+r\\n-j\end{bmatrix}q^{(k-1)\left(\binom{n}{2}-\binom{j}{2}\right)}F^{(k)}_{kj+r+1}(x;q)$$

$$+\sum_{j=1}^{n}(-1)^{n-j}\begin{bmatrix}n+r\\n-j-1\end{bmatrix}q^{(k-1)\left(\binom{n}{2}-\binom{j+1}{2}\right)+kj+r+1}F^{(k)}_{kj+r+1}(x;q)+(-1)^n\begin{bmatrix}n+r\\n-1\end{bmatrix}q^{(k-1)\binom{n}{2}+r+1}F^{(k)}_{r+1}(x;q)$$

$$=\sum_{j=0}^{n}(-1)^{n-j}\begin{bmatrix}n+r+1\\n-j\end{bmatrix}q^{(k-1)\left(\binom{n}{2}-\binom{j}{2}\right)}F^{(k)}_{kj+r+1}(x;q)=h(n,k,r+1).$$

For $r=0$ we get

$$h(n,k,0)=\sum_{j=0}^{n}(-1)^{n-j}\begin{bmatrix}n\\j\end{bmatrix}q^{(k-1)\left(\binom{n}{2}-\binom{j}{2}\right)}F^{(k)}_{kj}(x;q)$$

We want to show that $h(n,k,0)=x^n F^{(k)}_{(k-1)n}(x;q)$.

Let

$$a(n,i)=x^{n-i}F^{(k)}_{(k-1)n+i}(x;q).$$

Then we have

$$a(n,i) = x^{n-i}F^{(k)}_{(k-1)n+i}(x;q) = x^{n-i-1}\left(F^{(k)}_{(k-1)n+i+1}(x;q) - q^{(k-1)(n-1)+i}F^{(k)}_{(k-1)(n-1)+i}(x;q)\right)$$
$$= a(n,i+1) - q^{(k-1)(n-1)+i}a(n-1,i)$$

with $a(n,n) = F^{(k)}_{kn}(x;q)$.

This gives

$$a(n,i) = \sum_{j=0}^{n}(-1)^{n-j}\begin{bmatrix}n-i\\n-j\end{bmatrix}q^{\left(\binom{n}{2}-\binom{j}{2}\right)k-\left(i(n-j)+\binom{n}{2}-\binom{j}{2}\right)}F^{(k)}_{kj}(x;q) = x^{n-i}F^{(k)}_{(k-1)n+i}(x;q). \quad (53)$$

We must only verify that $a(n,n) = F^{(k)}_{kn}(x;q)$ and $a(n,i) = a(n,i+1) - q^{(k-1)(n-1)+i}a(n-1,i)$.

$$a(n,n) = \sum_{j=0}^{n}(-1)^{n-j}\begin{bmatrix}0\\n-j\end{bmatrix}q^{\left(\binom{n}{2}-\binom{j}{2}\right)k-\left(i(n-j)+\binom{n}{2}-\binom{j}{2}\right)}F^{(k)}_{kj}(x;q) = F^{(k)}_{kn}(x;q).$$

$$a(n,i+1) - a(n,i) = \sum_{j=0}^{n}(-1)^{n-j}\begin{bmatrix}n-i-1\\n-j\end{bmatrix}q^{\left(\binom{n}{2}-\binom{j}{2}\right)k-\left((i+1)(n-j)+\binom{n}{2}-\binom{j}{2}\right)}F^{(k)}_{kj}(x;q)$$

$$-\sum_{j=0}^{n}(-1)^{n-j}\begin{bmatrix}n-i\\n-j\end{bmatrix}q^{\left(\binom{n}{2}-\binom{j}{2}\right)k-\left(i(n-j)+\binom{n}{2}-\binom{j}{2}\right)}F^{(k)}_{kj}(x;q)$$

$$= \sum_{j=0}^{n}(-1)^{n-j-1}F^{(k)}_{kj}(x;q)\left(\begin{bmatrix}n-i\\n-j\end{bmatrix} - q^{n-j}\begin{bmatrix}n-i-1\\n-j\end{bmatrix}\right)q^{\left(\binom{n}{2}-\binom{j}{2}\right)k-\left(i(n-j)+\binom{n}{2}-\binom{j}{2}\right)}$$

$$= \sum_{j=0}^{n}(-1)^{n-j-1}F^{(k)}_{kj}(x;q)\begin{bmatrix}n-i-1\\n-1-j\end{bmatrix}q^{\left(\binom{n}{2}-\binom{j}{2}\right)k-\left(i(n-j)+\binom{n}{2}-\binom{j}{2}\right)} = q^{(k-1)(n-1)+i}a(n-1,i).$$

Therefore we get

$$a(n,0) = \sum_{j=0}^{n}(-1)^{n-j}\begin{bmatrix}n\\n-j\end{bmatrix}q^{\left(\binom{n}{2}-\binom{j}{2}\right)k-\left(\binom{n}{2}-\binom{j}{2}\right)}F^{(k)}_{kj}(x;q) = \sum_{j=0}^{n}(-1)^{n-j}\begin{bmatrix}n\\j\end{bmatrix}q^{(k-1)\left(\binom{n}{2}-\binom{j}{2}\right)}F^{(k)}_{kj}(x;q) = h(n,0,k).$$

**Proof of Proposition 7**

This follows as above from the identity $\sum_{k=0}^{\lfloor n/2 \rfloor}(-1)^k\begin{bmatrix}n\\k\end{bmatrix}Luc_{n-2k}(x) = x^n$. A proof can be found in [4], Theorem 3.1.

**Proof of Proposition 8**

We must show that

$$\sum_{j=1}^{n} \begin{bmatrix} n-1 \\ j-1 \end{bmatrix} x^j F_j(x,q;q) = \sum_{j=0}^{n-1} \begin{bmatrix} n \\ j+1 \end{bmatrix} x^{j+2} F_j(x,q;q).$$

Comparing the coefficients of $x^{2n-2k}$ it suffices to show that

$$\sum_{j=0}^{k} \begin{bmatrix} n \\ k-j \end{bmatrix} \begin{bmatrix} n-k-1 \\ j \end{bmatrix} q^{j^2} = \sum_{j=0}^{k} \begin{bmatrix} n-1 \\ k-j \end{bmatrix} \begin{bmatrix} n-k \\ j \end{bmatrix} q^{j^2}.$$

I want to thank Christian Krattenthaler for the following proof.

We write the sums in $q$ – hypergeometric form

$$\sum_{j=0}^{k} \begin{bmatrix} n \\ k-j \end{bmatrix} \begin{bmatrix} n-k-1 \\ j \end{bmatrix} q^{j^2} = {}_2\varphi_1 \begin{bmatrix} q^{-k}, q^{1+k-n} \\ q^{1-k+n} \end{bmatrix} ; q, q^n \end{bmatrix} \frac{(q^{1-k+n}; q)_k}{(q;q)_k}, \tag{54}$$

$$\sum_{j=0}^{k} \begin{bmatrix} n-1 \\ k-j \end{bmatrix} \begin{bmatrix} n-k \\ j \end{bmatrix} q^{j^2} = {}_2\varphi_1 \begin{bmatrix} q^{-k}, q^{k-n} \\ q^{-k+n} \end{bmatrix} ; q, q^n \end{bmatrix} \frac{(q^{-k+n}; q)_k}{(q;q)_k}. \tag{55}$$

Applying Heine's transformation ([6],(III.2))

$${}_2\varphi_1 \begin{bmatrix} a,b \\ c \end{bmatrix} ; q, z \end{bmatrix} = {}_2\varphi_1 \begin{bmatrix} \frac{abz}{c}, b \\ bz \end{bmatrix} ; q, \frac{c}{b} \end{bmatrix} \frac{\left(\frac{c}{b}, bz; q\right)_\infty}{(c, z; q)_\infty} \quad \text{to the right-hand side of}$$

$${}_2\varphi_1 \begin{bmatrix} q^{-k}, q^{1+k-n} \\ q^{1-k+n} \end{bmatrix} ; q, q^n \end{bmatrix} \frac{(q^{1-k+n}; q)_k}{(q;q)_k} = {}_2\varphi_1 \begin{bmatrix} q^{1+k-n}, q^{-k} \\ q^{1-k+n} \end{bmatrix} ; q, q^n \end{bmatrix} \frac{(q^{1-k+n}; q)_k}{(q;q)_k}$$

we get

$${}_2\varphi_1 \begin{bmatrix} q^{k-n}, q^{-k} \\ q^{-k+n} \end{bmatrix} ; q, q^{n+1} \end{bmatrix} \frac{(q^{1-k+n}; q)_k}{(q;q)_k} \frac{(q^{n+1}; q)_\infty (q^{n-k}; q)_\infty}{(q^{n-k+1}; q)_\infty (q^n; q)_\infty} = {}_2\varphi_1 \begin{bmatrix} q^{k-n}, q^{-k} \\ q^{-k+n} \end{bmatrix} ; q, q^{n+1} \end{bmatrix} \frac{(q^{1-k+n}; q)_k}{(q;q)_k} \frac{1-q^{n-k}}{1-q^n}$$

$$= {}_2\varphi_1 \begin{bmatrix} q^{k-n}, q^{-k} \\ q^{-k+n} \end{bmatrix} ; q, q^{n+1} \end{bmatrix} \frac{(q^{n-k}; q)_k}{(q;q)_k}.$$

Thus (54) = (55).